\documentclass[11pt]{article}        

\usepackage{amsmath,amssymb,amsthm,hyperref,theoremref,titlesec,lipsum}

\providecommand{\keywords}[1]{\textbf{\textbf{Keywords}:} #1}
\providecommand{\subjclass}[1]{\textbf{\textbf{Mathematics subject classification(2010)}:} #1}
\newcommand\mycom[2]{\genfrac{}{}{0pt}{}{#1}{#2}}
\theoremstyle{plain}
\newtheorem{thm}{Theorem}[section]

\newtheorem{rem}{Remark}[section]
\newtheorem{defn}{Definition}[section]
\numberwithin{equation}{section}

\newcommand{\Q}{\mathbb Q}

\newcommand{\R}{\mathbb R}

\renewcommand{\P}{\mathbb P}
\def\f{{\mathfrak{f}}}
\title{Simultaneous sign change and equidistribution of signs of Fourier coefficients of two cusp forms}

\author{Mohammed Amin Amri \footnote{ACSA Laboratory, Department of Mathematics, Faculty of Sciences, Mohammed First University, Oujda, Morocco amri.amine.mohammed@gmail.com} }
 \date{November 6, 2017}

\begin{document}

\maketitle

\begin{abstract}
We study the simultaneous sign change of Fourier coefficients of a pair of distinct normalized newforms of integral weight supported on primes power indices, we also prove some equidistribution results. Finally, we consider an analogous question for Fourier coefficients of a pair of half-integral weight Hecke eigenforms.
\end{abstract}

\subjclass{11F03, 11F30, 11F37}

\keywords{Simultaneous sign changes, Fourier coefficients, Sato-Tate conjecture}

\section{Introduction and statements of results}
Throughout the paper we denote by $\P$ the set of all prime numbers, for a set $S\subset \P$ we denote by $\delta(S)$ its natural density. We write $S_k^{\mathrm{new}}(N)$ for the space of newforms of even integral weight $k$ with trivial Dirichlet character over the group $\Gamma_0(N)$. The letter $\mathcal{H}$ stands for the upper half-plane, for $z\in\mathcal{H}$ we set $q:=e^{2\pi iz}$.

The study of sign changes of Fourier coefficients of a single cusp form has been the focus of much recent study (cf. \cite{amri1,amri3,bruin,kohnen10,inam2013,lau2016,Kohnen2013}) due to their various number theoretic applications (see for instance \cite{fung}) and his long history which goes back at least to Siegel \cite{siegel}. The question of simultaneous sign change of Fourier coefficients of modular forms was first considered by Kohnen and Sengupta in \cite{Kohnen2009}. Indeed, they proved that for given two cusp forms of different weights and the same level with real algebraic Fourier coefficients, there exists a Galois automorphism $\sigma$ such that $f_1^{\sigma}$ and $f_2^{\sigma}$ have infinitely many Fourier coefficients of opposite sign using analytic properties of the Rankin--Selberg zeta function, Landau's theorem on Dirichlet series with non-negative coefficients, and the bounded denominator principle. Recently, this result was extended to two cusp forms with arbitrary real Fourier coefficients in \cite{Gun2015} by Gun, Kohnen, and Rath. Moreover, they show for any pair $(f_1,f_2)$ of distinct newforms with trivial Dirichlet character that the sequence $\{\lambda_1(p^{\nu})\lambda_2(p^{\nu})\}_{\nu\in\mathbb{N}}$ changes signs infinitely often for an infinite set of primes $p$, where $\lambda_i(n)$ is the $n$-th Fourier coefficient of $f_i$. 

The primary result of this paper is to calculate the proportion of integers for which the product $\lambda_1(p^{\nu})\lambda_2(p^{\nu})$ has the same sign. More precisely, we shall show the following equidistribution result.
\begin{thm}\thlabel{thm1}
Let 
$$
f_i(z)=\sum_{n\ge 1} \lambda_i(n)n^{(k_i-1)/2}q^n\in S_{k_i}^{\mathrm{new}}(N_i)\quad i=1,2,
$$
be two distinct normalized newforms of level $N_1,N_2$ and weights $k_1,k_2\ge2$ respectively. For any prime $p\nmid N_i$, define $\theta_i(p)\in [0,\pi]$ by the relation $\lambda_i(p)=2\cos \theta_i(p)$ with $i=1,2$. Suppose that $1,\frac{\theta_1(p)}{2\pi},\frac{\theta_2(p)}{2\pi}$ are linearly independent over $\mathbb{Q}$, then we have
$$
\lim_{x\to\infty}\dfrac{\#\{v\le x : \lambda_1(p^\nu)\lambda_2(p^\nu)\lessgtr 0\}}{x}=\frac{1}{2}.
$$
\end{thm}
The strategy of the proof is inspired by the one used by Kohnen et al., in \cite[Proof of Theorem 3, case (iv)]{Kohnen2013} and involve essentially a two dimensional variant of Weyl's equidistribution theorem (see \thref{Weyl} for a precise statement). 
\begin{rem}
It should be noted that if one of the numbers $\frac{\theta_1(p)}{2\pi},\frac{\theta_2(p)}{2\pi}$ is rational and the other is irrational, then we have
$$
\lim_{x\to\infty}\dfrac{\#\{v\le x : \lambda_1(p^\nu)\lambda_2(p^\nu)\lessgtr 0\}}{\#\{v\le x : \lambda_1(p^\nu)\lambda_2(p^\nu)\ne 0\}}=\frac{1}{2},
$$ 
which may be proved in much the same way used by the author in \cite[Proof of Theorem 4]{amri2}.
\end{rem}
Our next purpose is to deal with the sign changes of the sequence $\{\lambda_1(p^\nu)\lambda_2(p^\nu)\}_{p\in\P}$. More precisely, we will show the following.
\begin{thm}\thlabel{thm2}
Suppose that $k_1,k_2\ge2$ and $N_1,N_2$ are distinct natural numbers. Suppose that
\[
f_i(z):=\sum_{n\ge1}\lambda_i(n)n^{(k_i-1)/2}q^n\in S_{k_i}^{\mathrm{new}}(N_i) \quad  i=1,2,
\]
are normalized newforms without complex multiplication which are not twists of each other. Let $\nu$ be an odd integer, we let $\P_{<0}(\nu)$ denote the set 
$$
\{p \in\P : p\nmid N,\;\lambda_1(p^\nu)\lambda_2(p^\nu)<0\},
$$
and similarly $\P_{\le 0}(\nu)$, $\P_{>0}(\nu)$, $\P_{\ge 0}(\nu)$, $\P_{=0}(\nu)$, where $N=N_1N_2$. Then the sequence $\{\lambda_1(p^\nu)$ $\lambda_2(p^\nu)\}_{p\in\P}$ changes signs infinitely often. Moreover, the sets
$$\P_{<0}(\nu),\;\P_{>0}(\nu),\; \P_{\le 0}(\nu),\;\P_{\ge 0}(\nu),$$
have natural density equal to $\frac{1}{2}$ and $\delta(\P_{=0}(\nu))=0$.
\end{thm}
Next, we investigate the simultaneous sign change of Fourier coefficients of two half-integral weight Hecke eigenforms. In order to state our result, we shall need to introduce some notations. Let $k,N$ be natural numbers, fix a Dirichlet character  $\chi$ modulo $4N$. We shall denote by $S_{k+1/2}(4N,\chi)$ the space of cusp forms of weight $k+1/2$ over the congruence subgroup $\Gamma_0(4N)$ with character $\chi$. When $k=1$, we shall work only with the orthogonal complement with respect to the Petersson scalar product of the subspace generated by the unary theta functions. In this set-up we have.
\begin{thm}\thlabel{thm3}
Suppose that $k_1, k_2$ are distinct natural numbers, $N_1,N_2$ are odd square-free natural numbers, and $\chi_1, \chi_2$ are real characters modulo $4N_1, 4N_2$, respectively. Suppose that
$$
\f_i(z)=\sum_{n\ge1}a_i(n)q^n\in S_{k_i+1/2}(4N_i,\chi_i)\quad i=1,2,
$$
are Hecke eigenforms. Let $\nu$ be a positive odd integer, and $t$ be a square-free integer such that $a_1(t) a_2(t)\ne0$. Assume that the Shimura lift of $\f_1$ and $\f_2$ are not twists of each other. Define the set of primes 
$$
\P'_{<0}(\nu):=\{p \in\P : p\nmid N,\;a_1(tp^{2\nu})a_2(tp^{2\nu})<0\},
$$
where $N=2N_1N_2$, and similarly
$$\P'_{<0}(\nu),\;\;\P'_{>0}(\nu),\;\; \P'_{\le 0}(\nu),\;\; \P'_{\ge 0}(\nu),\;\; \P'_{=0}(\nu).$$
Then the sequence $\{a_1(tp^{2\nu}) a_2(tp^{2\nu})\}_{p\in\mathbb{P}}$ changes signs infinitely often. Moreover, the sets
$$\P'_{<0}(\nu),\;\P'_{>0}(\nu),\; \P'_{\le 0}(\nu),\; \P'_{\ge 0}(\nu),$$
have natural density equal to $\frac{1}{2}$ and $\delta(\P'_{=0}(\nu))=0$.
\end{thm}
The above theorem improves the result of Kumar in \cite[Theorem 4.2]{Kumar}. Moreover, we can prove that a cuspidal Hecke eigenform $\f_i$ is determined uniquely by the sign of the sequence $\{a_i(tp^{2\nu})\}_p$ up to a positive constant using the technique of our proof.
\section{Preliminary results}
In this section, we introduce some notations and give the main necessary tools and definitions. Let 
\[
f_i(z):=\sum_{n\ge1}\lambda_i(n)n^{(k_i-1)/2}q^n \in S_{k_i}^{\mathrm{new}}(N_i),\quad i=1,2,
\]
be two distinct newforms of level $N_1, N_2$ and weights $k_1,k_2\ge 2$, respectively. In view of Deligne's bound, we have
$$
|\lambda_i(p)|\le 2\quad  i=1, 2,
$$
where $p$ is any prime number not dividing $N_i$. Thus for any prime $p\nmid N_i$ we can write
\begin{equation}\label{eq1}
\lambda_i(p)=2\cos\theta_i(p) \quad i=1,2,
\end{equation} 
for a uniquely defined angle $\theta_i(p)\in [0,\pi]$. Let us define the $2$-product Sato-Tate measure.
\begin{defn}
The $2$-product Sato-Tate measure $\mu_{\mathrm{ST}}^{\otimes_2}$ is the probability measure on $[0,\pi]^2$ given by 
$$\mu^{\otimes_2}_{\mathrm{ST}}:=\frac{4}{\pi^2}\sin^2\theta_1 \sin^2\theta_2 d\theta_1 d\theta_2.$$
\end{defn}
We now state the now-proven pair-Sato-Tate conjecture \cite[Proposition 2.2]{Wong} which will play an important role in proving \thref{thm1} and \ref{thm3}.
\begin{thm}[Pair-Sato-Tate conjecture]\thlabel{con1}
Let $k_1,k_2\ge2$ and $N_1,N_2$ are distinct natural numbers, and let
\[
f_i(z):=\sum_{n\ge1}\lambda_i(n)n^{(k_i-1)/2}q^n\in S_{k_i}^{\mathrm{new}}(N_i) \quad  i=1,2,
\]
are normalized newforms without complex multiplication, such that neither is a quadratic twist of the other. For a prime $p$, define $\theta_i(p)\in[0,\pi]$ to be as in \eqref{eq1}. Then the sequence $(\theta_1(p),\theta_2(p))$  is uniformly distributed in $[0,\pi]^2$ with respect to the $2$-product Sato-Tate measure $\mu_{\text{ST}}^{\otimes_2}$ as $p$ ranges over primes not dividing $N=N_1N_2$. In particular, for any two sub-interval $I_1\subset[0,\pi]$ $I_2\subset[0,\pi]$, we have 
$$
\lim_{x\to \infty}\dfrac{\#\{p\le x : p\nmid N,\;(\theta_1(p),\theta_2(p))\in I_1\times I_2 \}}{\#\{p\in\P  :  p\le x\}}=\mu^{\otimes_2}_{\mathrm{ST}}(I_1\times I_2)= \frac{4}{\pi^2}\iint_{I_1\times I_2}\sin^2\theta_1 \sin^2\theta_2 \,d\theta_1 \,d\theta_2.
$$
\end{thm}

In order to state a two dimensional variant of Weyl's equidistribution theorem which will be used to prove \thref{thm2}, we shall need to introduce notations and definitions following \cite{kuipers}. Let $\mathbf{a}=(a_1, a_2)$ and $\mathbf{b}=(b_1,b_2)$ be two vectors with real components. We say that $\mathbf{a}\le \mathbf{b}$ if $a_i\le b_i$ for $i=1,2$. The set of points $\textbf{x}\in\R^2$ such that $\textbf{a}\le \textbf{x} \le \textbf{b}$ will be denoted by $[\mathbf{a}, \mathbf{b}]$. We denote by $[\mathbf{0},\mathbf{1}]$, the two dimensional unit cube where $\mathbf{0}= (0,0)$ and $\mathbf{1}= (1,1)$. The integral part of $\mathbf{x} = (x_1,x_2)$ is $[\mathbf{x}] = ([x_1],[x_2])$ and the fractional part of $\mathbf{x}$ is $\{\mathbf{x}\} = (\{x_1\},\{x_2\})$. 
\begin{defn}
A sequence $(\mathbf{x}_n)_{n\in\mathbb{N}}$ of vectors in $\R^2$, is said to be uniformly distributed $\pmod 1$
in $\R^2$ if
$$
\lim_{N\to\infty}\dfrac{\#\{n\le N : \{\mathbf{x}_n\}\in [\mathbf{a}, \mathbf{b}]\}}{N}=(b_1-a_1)(b_2-a_2),
$$
for all sub-intervals $[\mathbf{a}, \mathbf{b}]\subseteq[\mathbf{0},\mathbf{1}]$.
\end{defn}
At this point we state the following theorem which is a straightforward consequence of \cite[Theorem 6.3, pp.48]{kuipers}.
\begin{thm}\thlabel{Weyl}
Let $\mathbf{\theta}=(\theta_1,\theta_2)\in\mathbb{R}^2$ be a vector so that the real numbers $1,\theta_1,\theta_2$ are linearly independent over $\Q$, then the sequence $(n\mathbf{\theta})=(n\theta_1,n\theta_2)$ is uniformly distributed $\pmod 1$ in $\R^2$.
\end{thm}

\section{Proofs}

\subsection{Proof of \texorpdfstring{\thref{thm1}}{Theorem 1.1}}
Since $f_i$ is a newform, then in particular it is a Hecke eigenform. Consequently its $p^\nu$-th eigenvalue is expressible by the following elementary trigonometric identity
$$
\lambda_i(p^\nu)=\dfrac{\sin((\nu+1)\theta_i(p))}{\sin\theta_i(p)}\quad\quad i=1,2,
$$
where $\theta_i(p)\in (0,\pi)$ be as in \eqref{eq1}, and obvious interpretation in the limiting cases $\theta_i(p)=0,\pi$, which can happen for only finitely many primes.

We set 
$\theta(p):=(\theta_1(p),\theta_2(p))$, $\nu\theta(p)=(\nu\theta_1(p),\nu\theta_2(p))$,
 and  $\sin(\nu\theta(p)):=(\sin(\nu\theta_1(p)),\sin(\nu\theta_2(p)))$. Writing 
$$
\nu\theta(p)=2\pi[\nu\theta(p)/2\pi]+2\pi\{\nu\theta(p)/2\pi\},
$$
it follows that
$$
\sin(\nu\theta(p))=\sin(2\pi\{\nu\theta(p)/2\pi\}).
$$
Thus for any sub-interval $[\mathbf{a},\mathbf{b}]\subset [\mathbf{-1},\mathbf{1}]$ one has
\begin{eqnarray*}
\sum_{\mycom{2\le \nu\le x+1}{\sin(\nu\theta(p))\in [\mathbf{a},\mathbf{b}]}}1 &=& \sum_{\mycom{2\le \nu\le x+1}{\sin(2\pi\{\nu\theta(p)/2\pi\}))\in [\mathbf{a},\mathbf{b}]}}1,\\
 &=&4\sum_{\mycom{2\le \nu\le x+1}{2\pi\{\nu\theta(p)/2\pi\}\in [\arcsin(\mathbf{a}),\arcsin(\mathbf{b})]}}1.
\end{eqnarray*}

By our assumption the numbers $1,\frac{\theta_1(p)}{2\pi},\frac{\theta_2(p)}{2\pi}$ are rationally independent. It follows from \thref{Weyl} that the sequence $\{\nu\theta(p)/2\pi\}$ is uniformly distributed $\pmod 1$ in $\R^2$. Consequently
\begin{equation}\label{eq3}
\sum_{\mycom{2\le \nu\le x+1}{\sin(\nu\theta(p))\in [\mathbf{a},\mathbf{b}]}}1\sim \dfrac{(\arcsin(b_1)-\arcsin(a_1))(\arcsin(b_2)-\arcsin(a_2))}{\pi^2} x,
\end{equation}
as $x\to\infty$. The claimed result now follows from \eqref{eq3} by taking into account that 
$$
\#\{\nu\le x : \lambda_1(p^\nu)\lambda_2(p^\nu)>0\}=\sum_{\mycom{2\le \nu\le x+1}{\sin(\nu\theta(p))\in [\mathbf{0},\mathbf{1}]}}1+\sum_{\mycom{2\le \nu\le x+1}{\sin(\nu\theta(p))\in [\mathbf{-1},\mathbf{0}]}}1,
$$
and
$$
\#\{\nu\le x : \lambda_1(p^\nu)\lambda_2(p^\nu)<0\}=\sum_{\mycom{2\le \nu\le x+1}{\sin(\nu\theta(p))\in [\mathbf{a},\mathbf{b}]}}1+\sum_{\mycom{2\le \nu\le x+1}{\sin(\nu\theta(p))\in [\mathbf{a}',\mathbf{b}']}}1,
$$
where $\mathbf{a}=(-1,0)$, $\mathbf{b}=(0,1)$, $\mathbf{a}'=(0,-1)$ and $\mathbf{b}'=(1,0)$. 
\subsection{Proof of \texorpdfstring{\thref{thm2}}{Theorem 1.2}}
As in the foregoing proof, we have
$$
\lambda_i(p^\nu)=\dfrac{\sin((\nu+1)\theta_i(p))}{\sin\theta_i(p)}\quad\quad i=1,2,
$$
where $\theta_i(p)\in[0,\pi]$. Since the set $\{p\in\mathbb{P} : \theta_i(p)=0,\pi\}$ is finite, we may assume that $\theta_i(p)\in(0,\pi)$. It follows that the sign of the sequence $\lambda_i(p^\nu)$ is the same as the sign of $\sin((\nu+1)\theta_i(p))$. Therefore
$$
\lambda_i(p^\nu)>0\quad\text{if and only if}\quad \theta_i(p)\in A_{>0}:=\bigcup_{j=1}^{\frac{\nu+1}{2}}\left(\dfrac{2(j-1)\pi}{\nu+1},\dfrac{(2j-1)\pi}{\nu+1}\right),
$$
and 
$$
\lambda_i(p^\nu)<0\quad\text{if and only if}\quad \theta_i(p)\in A_{<0}:=\bigcup_{j=1}^{\frac{\nu+1}{2}}\left(\dfrac{(2j-1)\pi}{\nu+1},\dfrac{2j\pi}{\nu+1}\right).
$$

By \thref{con1} the sequence $\{(\theta_1(p),\theta_2(p))\}_p$ is uniformly distributed in $[0,\pi]^2$ when $p$ ranges over primes numbers $p\nmid N$. Therefore, there exists infinitely many primes $p$ such that 
$$
(\theta_1(p),\theta_2(p))\in A_{>0}\times A_{>0},
$$
and there exists infinitely many primes $p$ such that 
$$
(\theta_1(p),\theta_2(p))\in  A_{>0}\times A_{<0}.
$$
Thus the sequence $\{\lambda_1(p^\nu)\lambda_2(p^\nu))\}_{p\in\P}$ changes sign infinitely often.
On the other hand, we have 
$$
\P_{>0}(\nu)=\{p\in \P : p\nmid N,\;(\theta_1(p),\theta_2(p))\in A_{>0}\times A_{>0} \} \sqcup\{p\in \P : p\nmid N,\; (\theta_1(p), \theta_2(p))\in A_{<0}\times A_{<0} \},
$$
and
$$
\P_{<0}(\nu)=\{p\in \P : p\nmid N,\;(\theta_1(p),\theta_2(p))\in A_{>0}\times A_{<0}\}\sqcup  \{p\in \P : p\nmid N,\;(\theta_1(p),\theta_2(p))\in A_{<0}\times A_{>0} \}.
$$
From \cite{Meher} we have 
$$
\mu_{\text{ST}}^{\otimes_2}(A_{>0}\times A_{<0})=\mu_{\text{ST}}^{\otimes_2}(A_{<0}\times A_{<0})=\frac{1}{4}.
$$
Therefore, by \thref{con1} we get
$$
\delta(\mathbb{P}_{>0}(\nu))=\delta(\mathbb{P}_{<0}(\nu))=\frac{1}{2}.
$$
An argument similar to the above shows that  $\delta(\mathbb{P}_{\ge 0}(\nu))=\delta(\mathbb{P}_{\le0}(\nu))=\frac{1}{2}$ and a fortiori we have $\delta(\mathbb{P}_{=0}(\nu))=0$.
\subsection{Proof of \texorpdfstring{\thref{thm3}}{Theorem 1.3}}
We will start off by recalling some generalities about Shimura correspondence from \cite{shimura73,niwa75}. The Shimura correspondence lifts $\f_i$ to a cusp form $f_{t,i}$ of weight $2k_i$ for the group $\Gamma_0(2N_i)$ with character $\chi_i^2$. Let 
$$
f_{t,i}(z)=\sum_{n\ge 1} A_{t,i}(n)q^n\quad (z\in\mathcal{H}),
$$
be its expansion at $\infty$. According to \cite{shimura73}, the $n$-th Fourier coefficient of $f_{t,i}$ is given by 
\begin{equation}\label{eq4}
A_{t,i}(n)=\sum_{d|n}\chi_{t,N_i}(d)d^{k_i-1}a_i\left(\frac{n^2}{d^2}t\right),
\end{equation}
where $\chi_{t,N_i}$ denotes the character $\chi_{t,N_i}(d):=\chi_i(d)\left(\frac{(-1)^{k_i}N_i^{2}t}{d}\right)$. Since $\f_i$ is a Hecke eigenform, then, so is the Shimura lift. Indeed, we have $f_{t,i}=a_i(t)f_i$ where $f_i$ is a normalized Hecke eigenform, write 
$$
f_i(z)=\sum_{n\ge 1}\lambda_i(n)n^{k_i-1/2}q^n \quad (z\in\mathcal{H}),
$$
for its Fourier expansion at $\infty$. We shall assume that $a_i(t)=1$, in the general case we may apply the proof to $\frac{\f_i}{a_i(t)}$. Since $2N_i$ is square-free, it follows that $f_i$ is a normalized Hecke eigenform without complex multiplication independent of $t$. 

Applying the M\"obius inversion formula to \eqref{eq4}, we derive that
$$
a_i(tn^2)=\sum_{d |n} \mu(d)\chi_{t,N_i}(d)d^{k_i-1} A_{t,i}\left(\frac{n}{d}\right).
$$
The above equality specialises to 
\begin{equation}\label{eq5}
\dfrac{a_i(tp^{2\nu})}{p^{\nu(k_i-1/2)}}=\lambda_i(p^\nu)-\frac{\chi_{t,N_i}(p)}{\sqrt{p}}\lambda_i(p^{\nu-1}),
\end{equation}
by taking $n=p^\nu$ and normalizing by $p^{\nu(k_i-1/2)}$. Now, rewrite \eqref{eq5} into
$$
\dfrac{a_i(tp^{2\nu})}{p^{\nu(k_i-1/2)}}=\frac{\sin((\nu+1)\theta_i(p))}{\sin\theta_i(p)}-\frac{\chi_{t,N_i}(p)}{\sqrt{p}}\frac{\sin(\nu\theta_i(p))}{\sin\theta_i(p)},
$$
where $\theta_i(p)\in [0,\pi]$. As the set $\{p\in\mathbb{P} : \theta_i(p)=0,\pi\}$ is finite (see \cite[Remark 2]{Kohnen2013}). We may assume that $\theta_i(p)\in (0,\pi)$. Consequently, we have
\begin{equation}\label{eq6}
a_i(tp^{2\nu})>0 \Longleftrightarrow \sin((\nu+1)\theta_i(p))>\frac{\chi_{t,Nçi}(p)}{\sqrt{p}}\sin(\nu\theta_i(p)),
\end{equation}
\begin{equation}\label{eq7}
a_i(tp^{2\nu})<0 \Longleftrightarrow \sin((\nu+1)\theta_i(p))<\frac{\chi_{t,Nçi}(p)}{\sqrt{p}}\sin(\nu\theta_i(p)).
\end{equation}
Let $\epsilon>0$, then for all $p>\frac{1}{\epsilon^2}$ we have $\left|\frac{\chi_{t,N_i}(p)}{\sqrt{p}}\sin(\nu\theta_i(p))\right|<\epsilon$. This together with \eqref{eq6} and \eqref{eq7} then implies
\begin{multline}\label{eq8}
\left\{p>\frac{1}{\epsilon^2} : p\nmid N, \sin((\nu+1)\theta_i(p))<-\epsilon, i=1,2\right\}\sqcup \\   
\left\{p>\frac{1}{\epsilon^2} : p\nmid N,\sin((\nu+ 1)\theta_i(p)) >\epsilon, i=1,2\right\}\subset\P'_{>0}(\nu).
\end{multline}
On the other hand we have
$$
\sin((\nu+1)\theta_i(p))>\epsilon\Longleftrightarrow \theta_i(p)\in I_{\epsilon},
$$
and
$$
\sin((\nu+1)\theta_i(p))<-\epsilon\Longleftrightarrow \theta_i(p)\in I'_{\epsilon}.
$$
Where 
$$
I_{\epsilon}:=\bigcup_{j=1}^{\frac{\nu+1}{2}}\left(\frac{(2j-2)\pi+\arcsin(\epsilon)}{\nu+1},\frac{(2j-1)\pi-\arcsin(\epsilon)}{\nu+1}\right),
$$
and
$$
I'_{\epsilon}:=\bigcup_{j=1}^{\frac{\nu+1}{2}}\left(\frac{(2j-1)\pi+\arcsin(\epsilon)}{\nu+1},\frac{2j\pi-\arcsin(\epsilon)}{\nu+1}\right).
$$
Thus \eqref{eq8} becomes
\begin{multline}\label{eq9}
\left\{p>\frac{1}{\epsilon^2} : p\nmid N,\; (\theta_1(p),\theta_2(p))\in I_{\epsilon}\times I_{\epsilon}\right\}\sqcup \\ \left\{p>\frac{1}{\epsilon^2} : p\nmid N,\; (\theta_1(p),\theta_2(p))\in I'_\epsilon\times I'_\epsilon\right\}\subset\P'_{>0}(\nu).
\end{multline}
From \thref{con1} we know that the sequence $(\theta_1(p),\theta_2(p))$ is equidistributed in $[0,\pi]^2$ with respect to the measure $\mu_{ST}^{\otimes_2}$. It follows that there are infinitely many primes $p$ such that $(\theta_1(p),\theta_2(p))\in I_{\epsilon}\times I_\epsilon$ and infinitely many primes $p$ such that $(\theta_1(p),\theta_2(p))\in I'_{\epsilon}\times I'_{\epsilon}$. Hence, the sets in \eqref{eq9} are infinite and consequently there exists infinitely many primes $p$ such that $a_1(tp^{2\nu})a_2(tp^{2\nu})>0$. In a similar way one can see that there exists infinitely many primes $p$ such that $a_1(tp^{2\nu})a_2(tp^{2\nu})<0$. 

Next, from \eqref{eq9} one has
$$
\pi_{>0}(x)+\pi\left(\frac{1}{\epsilon^2}\right)\ge S(I_\epsilon,I_\epsilon)(x)+S(I'_\epsilon,I'_\epsilon)(x),
$$
where $S(I_\epsilon,I_\epsilon)(x):=\#\{p\le x : p\nmid N,\;(\theta_1(p),\theta_2(p))\in I_\epsilon\times I_\epsilon\}$, $S(I'_\epsilon,I'_\epsilon)(x):=\#\{p\le x : p\nmid N,\; (\theta_1(p),\theta_2(p))\in I'_\epsilon\times I'_\epsilon\}$, and $\pi_{>0}(x):=\{p\le x : p\nmid N,\; a_1(tp^{2\nu})a_2(tp^{2\nu})>0\}$. Dividing the above inequality by $\pi(x)$ (the prime-counting function), we obtain
\begin{equation}\label{eq10}
\frac{\pi_{>0}(x)}{\pi(x)}+\dfrac{\pi\left(\frac{1}{\epsilon^2}\right)}{\pi(x)}\geq \dfrac{S(I_\epsilon,I_\epsilon)(x)+S(I'_\epsilon,I'_\epsilon)(x)}{\pi(x)}.
\end{equation}
By \thref{con1} we have 
\begin{equation}\label{eq11}
 \lim_{x\to\infty}\dfrac{S(I_\epsilon,I_\epsilon)(x)}{\pi(x)}=\mu_{\mathrm{ST}}^{\otimes_2}\left(I_\epsilon\times I_\epsilon\right)
 \;\;\text{and}\;\;
\lim_{x\to\infty}\dfrac{S(I'_\epsilon, I'_\epsilon)(x)}{\pi(x)}=\mu_{\mathrm{ST}}^{\otimes_2}\left(I'_\epsilon\times I'_\epsilon\right).
 \end{equation}
Since $\pi\left(\frac{1}{\epsilon^2}\right)$, is finite, the term $\frac{\pi\left(\frac{1}{\epsilon^2}\right)}{\pi(x)}$ tends to zero as $x$ tends to infinity. Taking into account \eqref{eq11}, a passage to the limit in \eqref{eq10} implies that
\begin{equation}\label{eq12}
\liminf_{x\to\infty}\frac{\pi_{>0}(x)}{\pi(x)}\geq \mu_{\mathrm{ST}}^{\otimes_2}\left(I_\epsilon\times I_\epsilon\right)+\mu_{\mathrm{ST}}^{\otimes_2}\left(I'_\epsilon\times I'_\epsilon\right).
\end{equation}
Letting $\epsilon$ to zero in \eqref{eq12}, we find 
$$
\liminf_{x\to\infty}\frac{\pi_{>0}(x)}{\pi(x)}\geq \mu_{\mathrm{ST}}^{\otimes_2}\left(I\times I\right)+\mu_{\mathrm{ST}}^{\otimes_2}\left(I'\times I'\right),
$$
where 
$$
I:=\bigcup_{j=1}^{\frac{\nu+1}{2}}\left(\frac{(2j-2)\pi}{\nu+1},\frac{(2j-1)\pi}{\nu+1}\right)\;\;\text{and}\;\;
I':=\bigcup_{j=1}^{\frac{\nu+1}{2}}\left(\frac{(2j-1)\pi}{\nu+1},\frac{2j\pi}{\nu+1}\right).
$$
From \cite{Meher} we have 
$$
\mu_{\mathrm{ST}}^{\otimes_2}\left(I\times I\right)=\mu_{\mathrm{ST}}^{\otimes_2}\left(I'\times I'\right)=\frac{1}{4}.
$$
It follows that 
$$
\liminf_{x\to\infty}\frac{\pi_{>0}(x)}{\pi(x)}\geq \frac{1}{2}.
$$
The same arguments yield $\liminf\limits_{x\to\infty}\frac{\pi_{\le0}(x)}{\pi(x)}\geq \frac{1}{2}$. 

Finally, upon using $\pi_{\le0}(x)=\pi(x)-\pi_{>0}(x)$, we deduce that 
$$
\frac{1}{2}\le \liminf\limits_{x\to\infty}\frac{\pi_{>0}(x)}{\pi(x)}\le \limsup\limits_{x\to\infty}\frac{\pi_{>0}(x)}{\pi(x)}\le \frac{1}{2},
$$
and hence $\lim\limits_{x\to\infty}\frac{\pi_{>0}(x)}{\pi(x)}$ exists and is equal to $\frac{1}{2}$. The rest of the proof runs as before. Consequently we have $\delta(\P'_{=0})=0$.
\subsection*{Acknowledgments}
The author would like to thank the referees for their careful reading and for pointing out the reference of Wong.

\titleformat{\section}[display]
{\normalfont\fillast}
{\scshape section \oldstylenums{\thesection}}
{1ex minus .1ex}
{\normalfont\large\bfseries}

% ------------------------------------------------------------------------
\end{document}